\documentclass[11pt,reqno,a4paper]{amsart}
\usepackage{euscript}

\newtheorem{theorem}{Theorem}

\newtheorem{lemma}{Lemma}

\renewcommand{\epsilon}{\varepsilon}

\def\Id{\text{\rm Id}}
\def\cA{\EuScript{A}}

\def\Z{\mathbb{Z}}
\def\R{\mathbb{R}}
\allowdisplaybreaks[1]

\numberwithin{equation}{section}

\begin{document}

\title[Smooth Linearization with Nonuniform Dichotomy]
{Smooth Linearization of Nonautonomous Differential Equations with a Nonuniform Dichotomy
}


\maketitle

\author{Davor Dragi\v cevi\'c\,$^a$, \ \ \ Weinian Zhang\,$^b$,\ \ \  Wenmeng Zhang\,$^c$
\footnote{Corresponding author: wmzhang@cqnu.edu.cn}
\\   \\
$^a${\small Department of Mathematics, University of Rijeka}
\\
\indent
{\small 51 000 Rijeka, Croatia, ddragicevic@math.uniri.hr}
\\
$^b${\small School of Mathematics, Sichuan University}
\\
\indent
{\small Chengdu, Sichuan 610064, China, matzwn@126.com}
\\
$^c${\small School of Mathematical Sciences, Chongqing Normal University}
\\
\indent
{\small Chongqing, 401331, China, wmzhang@cqnu.edu.cn}
}

\begin{abstract}
In this paper we give a smooth linearization theorem for nonautonomous differential equations with a nonuniform strong exponential dichotomy. In terms of a discretized evolution operator with hyperbolic fixed point 0, we formulate its spectrum and then give a spectral bound condition for the linearization of such equations to be simultaneously differentiable at 0 and H\"older continuous near 0. Restricted to the autonomous case, our result is the first one that gives a rigorous proof for simultaneously differentiable and H\"older linearization of hyperbolic systems without any non-resonant conditions.
\end{abstract}

{\footnotesize \sc Keywords:}
Nonautonomous differential equation; nonuniform strong exponential dichotomy; smooth linearization; spectral bound; non-resonant condition.

{\footnotesize \sc MSN (2010):}
37C60;  
37D25.  

\parskip 0.4cm


\section{Introduction}
\label{sec.1}

Linearization,
which answers whether a dynamical system is locally conjugated to its linear part in the sense of $C^r$ ($r\ge 0$),
is one of the most fundamental and important problems in the theory of dynamical systems
and a powerful tool in discussion of qualitative properties.
Earlier works were made for analytical conjugation in the complex case.
Poincar\'e (\cite{Poin}) proved that an analytic diffeomorphism can be analytically conjugated to
its linear part near a fixed point if all eigenvalues of the linear part lie inside the unit circle $S^1$ (or outside $S^1$) and satisfy the nonresonant condition. Siegel (\cite{Sieg}), Brjuno (\cite{Brj})
and Yoccoz (\cite{Yoc}) made contributions to the case of eigenvalues on $S^1$, in which the small divisor problem is involved.
On the other hand, the most well-known result in the real case is the Hartman-Grobman Theorem (\cite{Hart}), which says that
$C^1$ diffeomorphisms in $\mathbb{R}^n$ can be $C^0$ linearized near the hyperbolic fixed points. Later this result was generalized to
Banach spaces by  Palis (\cite{Palis}) and Pugh (\cite{Pugh}).

In order to improve smoothness of the conjugacy in the Hartman-Grobman Theorem to preserve more dynamical properties such as the characteristic direction (i.e., a fixed direction
along which an orbit approaches the equilibrium point), the smoothness of invariant manifold and the convergence (or divergence) rate of iteration, efforts were made to obtain smooth linearization.
In 1950's Sternberg (\cite{Stern58}) proved that $C^k$ ($k\ge1$) diffeomorphisms can be $C^r$ linearized near the hyperbolic fixed points, where the integer $r$ depends on $k$ and the non-resonant condition.
In 1970's Belitskii (\cite{Bel73}) gave conditions on $C^k$ linearization for $C^{k,1}$ ($k\ge 1$) diffeomorphisms, which implies that $C^{1,1}$ diffeomorphisms can be $C^1$ linearized locally if the eigenvalues $\lambda_1,...,\lambda_n$ satisfy a non-resonant condition that
\begin{align}\label{NRC}
|\lambda_i|\cdot|\lambda_j|\ne |\lambda_\iota|
\end{align}
for all $\iota=1,...,n$ if $|\lambda_i|<1<|\lambda_j|$.
This result was partially generalized to infinite-dimensional spaces in \cite{El-HJM14,R-S-JDDE04,ZZJ}. Notice that in the contractive (or expansive) case (\ref{NRC}) holds automatically and therefore $C^1$ linearization can always be realized in $\mathbb{R}^n$ (\cite{Hart60}). More results on $C^1$ linearization of contractions (or expansions) in Banach spaces can be found in \cite{ElB,Newhouse,R-S-JDE04,R-S-JDDE06}.
Concerning the general hyperbolic case,
although it is very important, Belitskii's $C^1$ linearization theorem cannot be used sometimes because the non-resonant condition (\ref{NRC}) may not be satisfied. Notice that (\ref{NRC}) cannot be removed for $C^1$ linearization due to a counterexample given by Hartman (\cite{Hart60}).
Thus, without any non-resonant conditions, most of the attentions were devoted to differentiable or/and H\"older linearization. On the one hand, although H\"older linearization of hyperbolic systems without any non-resonant conditions was known to some authors (see e.g. \cite{Stri-JDE90}), its rigorous proof was first published in the paper \cite{BDV} in 2007. On the other hand, differentiable linearization at the hyperbolic fixed point $0$ was proved in \cite{GuyHassRay-DCDS03} for $C^\infty$ diffeomorphisms in $\mathbb{R}^n$. This result was generalized to Banach spaces under a much weaker smoothness condition of $C^1$ plus $\alpha$-H\"older continuity (at the fixed point 0) together with a spectral bound condition in \cite{ZLZ}. It is worth mentioning that van Strien proved in \cite{Stri-JDE90} that the linearization of $C^2$ diffeomorphisms can be simultaneously differentiable at 0 and H\"older continuous near 0. However, his proof was pointed out to be wrong in \cite{Ray-JDE98}.

In this paper, we show that van Strien's result is true and we further extend his result to nonautonomous differential equations of the form
\begin{equation}\label{943}
x'=A(t)x+f(t,x)
\end{equation}
with the associated linear differential equation
\begin{equation}\label{944}
x'=A(t)x,
\end{equation}
where $A\colon \R \to M_d$ is a continuous map,
$M_d$ denotes the space of linear operators in $\R^d$,
and $f\colon \R \times \R^d \to \R^d$ is also a continuous
map such that $f(t, \cdot):\R^d \to \R^d$ is $C^{1}$. The first nonautonomous version of the Hartman-Grobman Theorem
for equation (\ref{943}) was given by Palmer (\cite{Palmer}) in 1973
under the assumption of (uniform) exponential dichotomy. More precisely,
assuming that~\eqref{944}
 admits a (uniform) exponential dichotomy and under appropriate assumptions for $f$, he proved
the $C^0$ linearization in the sense that
there exists a continuous function $H \colon \R \times \R^d \to \R^d$ such that
\begin{description}
\item[(H1)]
 for each $t\in \R$, $H_t:=H(t, \cdot)\colon \R^d \to \R^d$ is a homeomorphism;
\item[(H2)]
if $t\mapsto x(t)$ is a solution of~\eqref{943}, then $t\mapsto H(t, x(t))$ is a solution of~\eqref{944};
\item[(H3)]
 if $t\mapsto x(t)$ is a solution of~\eqref{944}, then $t\mapsto L(t, x(t))$,
 where
$L(t, x):=H_t^{-1}(x)$ for $t\in \R$ and $x\in \R^d$,
 is a solution of~\eqref{943}.
\end{description}
\vspace{-0.1cm}
After 2000, for hyperbolic nonautonomous differential equations, a result on linearization of Poincar\'e type with generalized non-resonant conditions was given in \cite{Sie-JDE02} and a result on H\"older linearization without any non-resonant conditons was given in \cite{BV2007}.
On the other hand, the problem of differentiable and smooth linearization for hyperbolic nonautonomous systems was considered only recently.  To the best of our knowledge, the first results in this direction were obtained in~\cite{CR}, where the authors formulated sufficient conditions for differentiable linearization of  nonautonomous dynamics whose linear part is uniformly exponentially stable.
More recently, this approach was extended to the case when the linear part of dynamics is  nonuniformly exponentially stable in~\cite{CMR}.  Furthermore, in~\cite{CDS} the authors have established a Sternberg type theorem for linear differential equations that gives conditions for smooth linearization of nonautonomous differential equations
whose linear part admits a uniform exponential dichotomy. In the case when the linear part of dynamics exhibits nonuniform exponential dichotomy, in~\cite{DZZ} we gave conditions for smooth linearization but for the case of discrete time.  In the present paper we formulate the first result that deals with smooth linearization under nonuniform hyperbolicity for continuous-time dynamics.
More precisely, we will extend van Strien's result of simultaneously differentiable and H\"older linearization to nonautonomous differential equations with a nonuniform exponential dichotomy. We emphasize that even in the uniform case our result is not covered by~\cite{CDS}. In fact, for a given integer $\ell\ge 1$, \cite{CDS} required (generalized) non-resonant conditions up to a sufficient larger order $k>\ell$ (i.e., all non-resonant conditions from order $2$ to order $k$) to guarantee the $C^\ell$ linearization. In comparison, our result only requires a spectral bound
condition, which is (using the same terminology as in~\cite{CDS}) a  type of non-resonant condition of order 2 (see details given just below the formulation of Theorem \ref{thm-mr}),
to guarantee simultaneously differentiable and H\"older linearization. Furthermore, the difference between the main  result in~\cite{CDS} and ours is even more obvious in the autonomous case,
where \cite{CDS} still needs the non-resonant conditions up to order $k$, whereas we do not need any non-resonant conditions since our spectral bound
condition holds automatically in the  autonomous case (see Theorem \ref{thm-aut}).

The strategy developed in this paper consists of passing from the continuous time nonautonomous dynamics to a discrete time nonautonomous dynamics. As in our previous work~\cite{DZZ}, we then study the associated autonomous dynamics on a larger space and after obtaining relevant results, we carefully return back to the setting of  nonautonomous dynamics.

We recall that in~\cite{DZZ} we dealt with perturbations of a linear dynamics
with discrete time that admits a nonuniform strong dichotomy.
In the present paper, we make the same assumption. Namely, we consider the case when our linear dynamics with continuous time
admits a nonuniform strong exponential dichotomy.
Recall that, in a definition of the ``exponential dichotomy", contractive and expansive parts of the evolution family of a linear system have bounded growth from above and below, respectively. However, since one needs to use the inverse of a  linear system for smooth linearization problem, ``strong" is imposed to guarantee that the contractive and expansive parts of the evolution family further have bounded growth from below and above, respectively (see Section 2 for more details).  We emphasize that all known versions of the Hartman-Grobman theorem in the nonautonomous setting  yielding the higher regularity of conjugacies (see~\cite{BV2007,BV,CDS}) require that the linear part of  the dynamics admits a strong exponential dichotomy. Indeed, although this terminology was not used in~\cite{CDS}, one observes that the condition~\cite[(A2)]{CDS}  implies that the authors assume that the linear part of the dynamics admits a strong uniform exponential dichotomy.

The paper is organized as follows.
We  formulate the spectrum of linear equation (\ref{944}) in terms of a discretized evolution operator in Section 2. Then we present our main theorem in Section 3, where a spectral
bound condition is given for the linearization of equation (\ref{943}) to be simultaneously differentiable at 0 and H\"older continuous near 0. Section 4 is devoted to the proof of our main theorem.
Finally, we extend our results to infinite-dimensional spaces in Section 5, which was not considered in \cite{CDS}.


\section{Dichotomies and spectrum
}
\label{sec.2}

Let $T(t,s)$ be an evolution family of the linear system \eqref{944}.
Following  \cite{BV},
we say that~\eqref{944} admits a \emph{nonuniform strong exponential dichotomy}
if
\begin{description}
\item[(SNED1)]
there is a family $P(t)$, $t\in \R$, of projections on $\R^d$ such that
\begin{eqnarray}
T(t,s)P(s)=P(t)T(t,s)\qquad \mbox{for }t,s \in \R,
\label{945}
\end{eqnarray}

\item[(SNED2)] there exist  $M, \lambda, \bar{\lambda} >0$, $\lambda \le \bar{\lambda}$  and $\epsilon \ge 0$ such that
\begin{eqnarray}
\left\{\begin{array}{rl}
\lVert T(t,s)P(s)\rVert \le Me^{-\lambda (t-s)+\epsilon \lvert s\rvert} \quad &\text{for $t\ge s$,}
\\
\lVert T(t,s)(\Id-P(s))\rVert \le Me^{-\lambda (s-t)+\epsilon \lvert s\rvert} \quad &\text{for $t\le s$}
\end{array}\right.
\label{ED1}
\end{eqnarray}
and
\begin{equation}\label{nbg}
\lVert T(t,s)\rVert \le Me^{\bar{\lambda} \lvert t-s\rvert+\epsilon \lvert s\rvert} \quad \text{for $t,s \in \R$.}
\end{equation}
\end{description}
This notion of dichotomy, introduced by Barreira and Valls~\cite{BV}, looks similar to
the concept of a  well-known nonuniform exponential dichotomy,
but the difference is that,
besides requiring exponential contraction along stable direction forward in time
and the exponential contraction along unstable direction backward in time (see~(\ref{ED1})),
it requires the evolution family $T(t,s)$ to exhibit the so-called nonuniform bounded growth condition (see \eqref{nbg})
.

We say that~\eqref{944} admits a \emph{strong exponential dichotomy with respect to
a family of norms} $\lVert \cdot \rVert_t$, $t\in \R$, if
\begin{description}
\item[(SED'1)]  there exists a family $P(t)$, $t\in \R$, of projections on $\R^d$ satisfying~\eqref{945},

\item[(SED'2)]
there exist $D, \lambda. \bar{\lambda} >0$, $\lambda \le \bar{\lambda}$ such that for each $x\in \R^d$,
\begin{equation}\label{d2}
\left\{\begin{array}{rl}
\lVert T(t,s)P(s)x\rVert_t \le Me^{-\lambda (t-s)}\lVert x\rVert_s \quad &\text{for $t\ge s$,}
\\
\lVert T(t,s)(\Id-P(s))x\rVert_t \le Me^{-\lambda (s-t)}\lVert x\rVert_s \quad &\text{for $t\le  s$}
\end{array}\right.
\end{equation}
and
\begin{equation}\label{d1}
\lVert T(t,s)x\rVert_t \le Me^{\bar{\lambda} \lvert t-s\rvert} \lVert x\rVert_s \quad \text{for $t,s \in \R$.}
\end{equation}
\end{description}
The following lemma gives  a relationship between those two concepts of dichotomy. It  is essentially established
in the proof of~\cite[Theorem 3.9]{BDV3}.

\begin{lemma}
\label{1038}
The following assertions are equivalent:
\begin{enumerate}
\item \eqref{944} admits a nonuniform strong exponential dichotomy;
\item \eqref{944}  admits a strong exponential dichotomy with respect to a family of norms $\lVert \cdot \rVert_t$, $t\in \R$ with the property that there exist $C>0$ and $\epsilon \ge 0$ such that
\begin{equation}\label{ln}
\lVert x\rVert \le \lVert x\rVert_t \le Ce^{\epsilon \lvert t\rvert} \lVert x\rVert,
\  \  \  \forall x\in \R^d \mbox{ and } t\in \R.
\end{equation}
\end{enumerate}
\end{lemma}


Let
\begin{equation}\label{an}
A_n:=T(n+1, n) \quad \text{for $n\in \Z$,}
\end{equation}
be the discretization of the evolution operator $T(t,s)$.
If \eqref{944} admits a nonuniform strong exponential dichotomy, then by Lemma~\ref{1038}
and \eqref{d2}-\eqref{d1} we see that
\begin{align}\label{Amn}
\cA(m, n):=\begin{cases}
A_{m-1}\cdots A_n & \text{for $m>n$,}\\
\Id & \text{for $m=n$,} \\
A_m^{-1} \cdots A_{n-1}^{-1} & \text{for $m<n$}
\end{cases}
\end{align}
admits a strong exponential dichotomy, i.e.,
for each $x\in \R^d$
\begin{eqnarray*}
\left\{\begin{array}{rl}
\lVert \cA(m, n)P(n)x\rVert_m  \le Me^{-\lambda (m-n)}\lVert x\rVert_n \quad \text{for $m\ge n$,}
\\
\lVert \cA(m, n)(\Id-P(n))x\rVert_m \le Me^{-\lambda (n-m)} \lVert x\rVert_n \quad \text{for $m\le n$}
\end{array}\right.
\end{eqnarray*}
and
\begin{equation}\label{1048}
\lVert \cA(m, n)x\rVert_m \le Me^{\bar{\lambda} \lvert m-n\rvert}\lVert x\rVert_n  \quad \text{for $m, n\in \Z$,}
\end{equation}
where
$\lVert \cdot \rVert_n$, $n\in \Z$, is a sequence of norms such that
\begin{equation}\label{ln2}
\lVert x\rVert \le \lVert x\rVert_n \le Ce^{\epsilon \lvert n\rvert} \lVert x\rVert,
\ \ \ \forall x\in \R^d \mbox{ and } n \in \Z,
\end{equation}
by (\ref{ln}). Let
\[
Y_\infty:=\bigg{\{} \mathbf x=(x_n)_{n\in \Z} \subset \R^d: \sup_{n\in \Z}\lVert x_n\rVert_n <\infty \bigg{\}}.
\]
Then, $(Y_\infty, \lVert \cdot \rVert)$ is a Banach space
equipped with the norm
$\lVert \mathbf x\rVert:=\sup_{n\in \Z}\lVert x_n\rVert_n$.
Define a linear operator $\mathbb A\colon Y_\infty \to Y_\infty$ by
\begin{equation}\label{mathbbA}
(\mathbb A\mathbf x)_n=A_{n-1}x_{n-1}, \quad \mathbf x=(x_n)_{n\in \Z}\in Y_\infty, \ n\in \Z.
\end{equation}
It follows from~\eqref{1048} that $\mathbb A$ is a well defined and bounded linear operator. Furthermore, $\mathbb A$ is invertible and
\[
(\mathbb A^{-1} \mathbf x)_n=A_n^{-1}x_{n+1}, \quad \mathbf x=(x_n)_{n\in \Z}\in Y_\infty, \ n\in \Z.
\]
We recall the following result.

\begin{lemma}{\rm (}Theorem {\rm 1} in \cite{DZZ}{\rm )}
\label{ts}
Let \eqref{944} admit a nonuniform strong exponential dichotomy and let $\mathbb{A}$ be defined in
{\rm (\ref{mathbbA})}.
Then there exist constants
\[
0< a_1\le b_1 <a_2\le b_2 <\ldots <a_k\le b_k<1<a_{k+1}\le b_{k+1}<\ldots <a_r\le b_r,
\]
such that
\begin{align}\label{nxx}
\sigma(\mathbb A)=\bigcup_{i=1}^r \{z\in \mathbb C: a_i\le \lvert z\rvert \le b_i\},
\end{align}
where $\sigma (\mathbb A)$ denotes the spectrum of $\mathbb A$.
\end{lemma}

It is worthy noting that we can describe $\sigma (\mathbb A)$ solely in terms of $T(t,s)$.
For each $\mu \in \mathbb R\setminus \{0\}$,  we can define a new evolution family $T_\mu (t,s)$ by
\[
T_\mu(t,s)=\frac{1}{\mu^{t-s}}T(t,s) \quad \text{for $t, s\in \R$.}
\]
Let $\lVert \cdot \rVert_t$, $t\in \R$, be the family of norms given by Lemma~\ref{1038}.

\noindent {\bf Proposition 1.}
{\it
$
\sigma (\mathbb A) \cap \R
$
is the set of all $\mu \in \mathbb R\setminus \{0\}$ such that $T_\mu (t,s)$ does not
admits a strong exponential dichotomy
with respect to
$\lVert \cdot \rVert_t, ~t\in \R$.
}

{\it Proof}.
Let $A_n$ be defined by~\eqref{an} for $n\in \Z$.
Assume that $\mu \in \mathbb R\setminus \{0\}$ is such that $T_\mu (t,s)$ admits a strong exponential dichotomy with respect to $\lVert \cdot \rVert_t$, $t\in \R$. This trivially implies that the sequence $(\frac{1}{\mu}A_m)_{m\in \Z}$ admits  a strong exponential dichotomy with respect to the sequence of norms
$\lVert \cdot \rVert_m$, $m\in \Z$.  Hence, \cite[Lemma 2]{DZZ} implies that $\mu \notin \sigma (\mathbb A)$.

Conversely, suppose that $\mu \notin \sigma (\mathbb A)$. Then, \cite[Lemma 2]{DZZ} implies that the sequence $(\frac{1}{\mu} A_m)_{m\in \Z}$ admits a strong exponential dichotomy with respect to the sequence of norms $\lVert \cdot \rVert_m$, $m\in \Z$. Let $P(n)$, $n\in \Z$ be the associated projections. A simple computation show that
$T_\mu(t,s)$ admits a strong exponential dichotomy with respect to the family of norms $\lVert \cdot \rVert_t$, $t\in \R$ and projections $P(t)$, $t\in \R$ given by
\[
P(t)=T(t,n)P(n)T(n, t) \quad \text{for $t\in [n, n+1)$ and $n\in \Z$.}
\]
The proof is complete.\qquad$\Box$


\section{Simultaneously Differentiable and H\"older Linearization}


In order to consider the simultaneously differentiable and H\"older linearization, we need to assume that the linear equation
\eqref{944} admits a nonuniform strong exponential dichotomy and therefore the spectrum $\sigma(\mathbb{A})$ for Eq.~\eqref{944} has the decomposition given in Lemma~\ref{ts}.
Moreover, we further assume that
$f\colon \R \times \R^d \to \R^d$ in~\eqref{943} is continuous and that  $f(t,\cdot)\colon \R^d \to \R^d$ is $C^1$ such that $D_xf(t,x)$ is a jointly continuous function of $(t,x)$ and
\begin{description}
\item[(F1)]
$f(t, 0)=0$ for all $t\in \R$;

\item[(F2)]
$D_xf(t, 0)=0$ for all $t\in \R$;

\item[(F3)]
$\lVert D_x f(t, x)\rVert \le \eta e^{-3\epsilon \lvert t\rvert}$ for all $t\in \R$, where
$\eta >0$ is a constant;

\item[(F4)]
$\lVert D_xf(t,x)-D_xf(t, y)\rVert \le Be^{-4\epsilon \lvert t\rvert} \lVert x-y\rVert$
for all $t\in \R$, where
$B>0$ is a constant.
\end{description}
Then we have the following main theorem of this paper.

\begin{theorem}\label{thm-mr}
Let \eqref{944} admit a nonuniform strong exponential dichotomy and let $a_1,...,a_r$ and $b_1,...,b_r$ be given in Lemma~{\rm \ref{ts}}. Assume that the spectral bound condition
 \begin{equation}\label{DH-cond}
b_i/a_i<b_k^{-1}, \ \forall i=1, \ldots, k, \ b_j/a_j < a_{k+1}, ~~~\forall j=k+1, \ldots, r
\end{equation}
holds and that $\alpha\in \mathbb{R}$ is an arbitrarily given constant satisfying
\begin{align}\label{def-a}
    0<\alpha<\min\Big\{\frac{\ln a_{k+1}-\ln b_k}{\ln b_r},\frac{\ln a_{k+1}-\ln b_{k}}{\ln a_1^{-1}}\Big\}.
\end{align}
Furthermore, suppose that
$f$ satisfies~{\bf (F1)}-{\bf (F4)}
with a sufficiently small constant $\eta>0$ {\rm (}which tends to $0$ when $\alpha$ tends to its upper bound{\rm )} and a constant
$B>0$.
Then there exist neighborhoods $V_t:=\{u\in\R^d:\|u\|\le e^{-2\epsilon|t|}\tilde \rho\}$ with a small constant $\tilde \rho>0$ {\rm(}independent of $\alpha${\rm )} and maps $H, G \colon  \R \times \R^d\to \R^d$ such that
\vspace{-0.1cm}
\begin{itemize}
\item[{\bf(A1)}] $H(t,x)=x+e^{(3+\varrho)\epsilon \lvert t\rvert} o(\lVert x\rVert^{1+\varrho})$, $G(t,x)=x+e^{(3+\varrho)\epsilon \lvert t\rvert} o(\lVert x\rVert^{1+\varrho})$, where $\varrho\in (0,1)$ is a small constant;
\vspace{0.2cm}
\item[{\bf(A2)}] $\|H(t,x)-H(t,y)\|\le \tilde C e^{(2+\alpha)\epsilon|t|} \|x-y\|^\alpha$ and $\|G(t,x)-G(t,y)\|\le \tilde C e^{(2+\alpha)\epsilon|t|} \|x-y\|^\alpha$ for all $x,y\in V_t$, where $\tilde C>0$ is a constant independent of $\alpha$;
\vspace{0.2cm}
\item[{\bf(A3)}] $H(t, G(t, x))=x$ and $G(t, H(t,x))=x$ for each $t\in \R$ and $x\in \R^d$;
\vspace{0.2cm}
\item[{\bf(A4)}] if $t\mapsto x(t)$ is a solution of~\eqref{943}, then $t\mapsto H(t, x(t))$ is a solution of~\eqref{944};
\vspace{0.2cm}
\item[{\bf(A5)}] if $t\mapsto x(t)$ is a solution of~\eqref{944}, then $t\mapsto G(t, x(t))$ is a solution of~\eqref{943}.
\end{itemize}
\end{theorem}
\vspace{-0.3cm}
Before giving the proof of Theorem~\ref{thm-mr},  we would like to compare it with the main result in~\cite{CDS}. Firstly, let us assume that~\eqref{944} admits a  uniform strong exponential dichotomy, i.e. nonuniform strong exponential dichotomy with $\epsilon=0$. In this case, $\sigma(\mathbb A)$ is closely related to
the so-called Sacker-Sell spectrum~\cite{SS} (see also~\cite{CDS} and
references therein for more details), which is denoted by $\Sigma_{SS}(A)$. More precisely, we have
$$
\Sigma_{SS}(A)=\bigcup_{i=1}^{r}[\ln a_i,\ln b_i].
$$
Observe that  our spectral gap condition (\ref{DH-cond}) is equivalent to requiring
\begin{align*}
&[\ln a_i,\ln b_i]\cap \{[\ln a_i,\ln b_i]+[\ln a_k,\ln b_k]\}=\emptyset,\quad
\forall i=1,...,k,
\\
&[\ln a_j,\ln b_j]\cap \{[\ln a_j,\ln b_j]+[\ln a_{k+1},\ln b_{k+1}]\}=\emptyset,\quad
\forall j=k+1,...,r.
\end{align*}
Using the terminology from~\cite{CDS}, the above condition is the so-called
 non-resonant condition of order $2$,
which is weaker than the non-resonant conditions up to a sufficiently
larger order $k$  required in~\cite[Theorem 5]{CDS}. However, our smoothness of simultaneously differentiable and H\"older continuity for linearization is lower than  $C^\ell$-smoothness ($\ell\ge 1$)  obtained in~\cite[Theorem 5]{CDS}.

{\it Proof of Theorem \ref{thm-mr}.}
Since $f(t,x)$ is continuous in $t$ and continuously differentiable in $x$ such that $f(t,0)=0$ for all $t\in\mathbb{R}$, and satisfies a global Lipschitz condition in $x$ with a constant $\eta>0$, by \cite[Theorems 1.6 and 2.4]{ZZF}, the solution $x(t):=\phi(t,t_0;x_0)$ of~\eqref{943} with $x(t_0)=x_0$ exists for all $t\in (-\infty,\infty)$ and is $C^1$ in $x_0$. Then we may define
$C^1$ maps $f_n \colon \R^d \to \R^d$ by
\begin{equation}\label{xfn}
f_n(x):=\phi(n+1, n; x) -A_nx\quad \text{for $x\in \R^d$.}
\end{equation}
We claim that $f_n$ satisfies
\begin{align}\label{fn0}
f_n(0)=0,\quad Df_n(0)=0 \quad {\rm for}~n\in\mathbb{Z},
\end{align}
and
\begin{align}\label{Dfn0}
\|D f_{n}(x)\| \le \tilde \eta e^{-\epsilon |n|},\quad
\|Df_{n}(x)-Df_{n}(y)\| \le \tilde B e^{-\epsilon |n|}\|x-y\|
\end{align}
for $x, y\in \R^d$ and $m\in \Z$, where $\tilde \eta,\tilde B>0$ are constants. In fact, $f_n(0)=0$ is clear by the fact that $\phi(t,t_0;0)=0$ for all $t$ and $t_0$ since $0$ is a solution, as known from {\bf (F1)}. For the derivative, we note that the variation of parameter formula implies that
\[
f_n(x)=\int_n^{n+1}T(n+1, r)f(r, \phi(r, n;x))\, dr,
\]
and therefore it follows from assumption {\bf (F2)} that
\[
\begin{split}
Df_n(0) &=\int_n^{n+1}T(n+1, r)D_xf(r, \phi(r, n;0))D_x\phi(r, n;0)\, dr \\
&=\int_n^{n+1}T(n+1, r)D_xf(r, 0)D_x\phi(r, n;0)\, dr=0.
\end{split}
\]
Thus the claimed result (\ref{fn0}) is proved.

For (\ref{Dfn0}), we observe that
\[
D_x\phi(t, n; x)=T(t,n)+\int_n^t T(t,r)D_xf(r, \phi(r, n;x))D_x\phi(r, n;x)\, dr
\]
for $t\ge n$ and $x\in \R^d$.  Hence, it follows from~\eqref{nbg} and
assumption {\bf (F3)}
that
\[
\begin{split}
&\lVert D_x\phi(t, n; x) \rVert
\\
&\le  Me^{\bar{\lambda} (t-n)+\epsilon \lvert n\rvert} + \int_n^t Me^{\bar{\lambda} (t-r)+\epsilon \lvert r\rvert}\eta e^{-3\epsilon \lvert r\rvert} \lVert D_x\phi(r, n;x)\rVert\, dr
\\
&\le Me^{\bar{\lambda}} e^{\epsilon \lvert n\rvert} +M\eta e^{\bar{\lambda }+2\epsilon}\int_n^t e^{-2\epsilon \lvert n\rvert} \lVert D_x\phi(r, n;x)\rVert\, dr
\end{split}
\]
for $t\in [n, n+1]$ and $x\in \R^d$. Hence, by Gronwall's lemma we get
\[
\begin{split}
\lVert D_x\phi(t, n; x) \rVert &\le Me^{\bar{\lambda}} e^{\epsilon \lvert n\rvert} e^{M\eta e^{\bar{\lambda} +2\epsilon} \int_n^t e^{-2\epsilon \lvert n\rvert} \, dr}
\end{split}
\]
and therefore
\begin{equation}\label{910x}
\lVert D_x\phi(t, n; x) \rVert \le \tilde M e^{\epsilon \lvert n\rvert}
\end{equation}
for every $t\in [n, n+1]$ and $x\in \R^d$, where $\tilde M:= Me^{\bar{\lambda}}e^{M e^{\bar{\lambda} +2\epsilon}}$ ($\eta$ can be removed here since it is small, i.e., $\eta<1$).
On the other hand, note that
\begin{equation}\label{Dfm}
Df_n(x)=\int_n^{n+1}T(n+1, r)D_xf(r, \phi(r, n;x))D_x\phi(r, n;x)\, dr.
\end{equation}
Then, combining (\ref{910x}) with (\ref{Dfm}) we get
\[
\begin{split}
\lVert Df_n(x) \rVert &\le \int_n^{n+1}Me^{\bar{\lambda} (n+1-r)+\epsilon \lvert r\rvert} \eta e^{-3\epsilon \lvert r\rvert}\tilde M e^{\epsilon \lvert n\rvert} \, dr \\
&\le M \tilde M \eta e^{\bar{\lambda}+2\epsilon+1}e^{-\epsilon \lvert n+1\rvert}
\end{split}
\]
for each $n\in \Z$ and $x\in X$. Hence, the first inequality of (\ref{Dfn0}) holds with
\begin{equation}\label{eta}
\tilde \eta:=M \tilde M \eta e^{\bar{\lambda}+2\epsilon+1} >0.
\end{equation}

For the second inequality of (\ref{Dfn0}), we observe that
\[
\begin{split}
&\phi(r,n;x)-\phi(r,n;y)
\\
&=T(r,n)(x-y)
+\int_n^r T(r,s)(f(s, \phi(s, n; x))-f(s, \phi(s, n; y)))\, ds
\end{split}
\]
for $r\ge n$ and $x,y\in \R^d$.
Thus, it follows from~\eqref{nbg} and assumption {\bf (F3)}
that
\[
\begin{split}
&\lVert \phi(r,n;x)-\phi(r,n;y)\rVert
\\
&\le Me^{\bar{\lambda} \lvert r-n\rvert+\epsilon \lvert n\rvert}\lVert x-y\rVert +\int_n^r M\eta e^{\bar{\lambda} \lvert r-s\rvert -2\epsilon \lvert s\rvert}\lVert \phi(s, n; x)-\phi(s, n; y)\rVert\, ds
\\
&\le Me^{\bar{\lambda} +\epsilon \lvert n\rvert}\lVert x-y\rVert +\int_n^r Me^{\bar{\lambda}} \eta  e^{-2\epsilon \lvert n\rvert +2\epsilon}\lVert \phi(s, n; x)-\phi(s, n; y)\rVert\, ds
\end{split}
\]
for every $r\in [n, n+1]$ and $x, y\in \R^d$. Then, Gronwall's lemma implies that there exists $a>0$ such that
\begin{equation}\label{912}
\lVert \phi(r,n;x)-\phi(r,n;y)\rVert \le ae^{\epsilon \lvert n\rvert}\lVert x-y\rVert
\end{equation}
for $n\in \Z$, $r\in [n, n+1]$ and $x, y\in \R^d$.
On the other hand, we have
\[
\begin{split}
&D_x\phi(t, n;x)-D_x\phi(t, n; y)
\\
&=\int_n^t T(t,r)D_xf(r, \phi(r, n;x))D_x\phi(r, n;x)\, dr \\
&\phantom{=}-\int_n^t T(t,r)D_xf(r, \phi(r, n;y))D_x\phi(r, n;y)\, dr \\
&=\int_n^t T(t,r)D_xf(r, \phi(r, n;x))(D_x\phi(r, n;x)-D_x\phi(r, n;y))\, dr \\
&\phantom{=}+\int_n^t T(t,r)(D_xf(r, \phi(r, n;x))-D_xf(r, \phi(r, n;y)))D_x\phi(r, n;y)\, dr.
\end{split}
\]
Hence, it follows from~\eqref{nbg}, assumptions {\bf (F3)} and
{\bf (F4)},
\eqref{910x} and~\eqref{912} that
\[
\begin{split}
&\lVert D_x\phi(t, n;x)-D_x\phi(t, n; y)  \rVert
\\
&\le \int_n^t Me^{\bar{\lambda} +\epsilon \lvert r\rvert}Be^{-3\epsilon \lvert r\rvert}\lVert \phi(r, n;x)-\phi(r, n;y)\rVert \tilde M e^{\epsilon \lvert n\rvert}\, dr\\
&\phantom{\le}+\int_n^t Me^{\bar{\lambda} +\epsilon \lvert r\rvert} \eta e^{-3\epsilon \lvert r\rvert}\lVert D_x\phi(r, n;x)-D_x\phi(r, n;y)\rVert \, dr \\
&\le \int_n^t Me^{\bar{\lambda} +\epsilon \lvert r\rvert}Be^{-3\epsilon \lvert r\rvert}ae^{\epsilon \lvert n\rvert}\lVert x-y\rVert  \tilde M e^{\epsilon \lvert n\rvert}\, dr\\
&\phantom{\le}+\int_n^t Me^{\bar{\lambda} +\epsilon \lvert r\rvert} \eta e^{-3\epsilon \lvert r\rvert}\lVert D_x\phi(r, n;x)-D_x\phi(r, n;y)\rVert \, dr
\end{split}
\]
for $t\in [n, n+1]$ and $x, y\in \R^d$. By Gronwall's inequality again, one can conclude that there exists $d>0$ such that
\begin{equation}\label{Df-Df}
\lVert D_x\phi(t, n;x)-D_x\phi(t, n; y)  \rVert  \le d\lVert x-y\rVert
\end{equation}
for $n\in \Z$, $t\in [n, n+1]$ and $x, y\in \R^d$.

Now we are ready to estimate the term $Df_{n}(x)-Df_{n}(y)$.
Since
\[
\begin{split}
&Df_{n}(x)-Df_{n}(y)
\\
&=\int_n^{n+1}T(n+1, r)D_xf(r, \phi(r, n;x))D_x\phi(r, n;x)\, dr \\
&\phantom{=}-\int_n^{n+1}T(n+1, r)D_xf(r, \phi(r, n;y))D_x\phi(r, n;y)\, dr \\
&=\int_n^{n+1}T(n+1, r)D_xf(r, \phi(r, n;x))(D_x\phi(r, n;x)-D_x\phi(r, n;y))\, dr \\
&\phantom{=}+\int_n^{n\!+\!1}T(n+1, r)(D_xf(r, \phi(r, n;x))
\!-\!D_xf(r, \phi(r, n;y)))D_x\phi(r, n;y)\, dr, \\
\end{split}
\]
we obtain from {\bf (F3)}-{\bf (F4)}, (\ref{910x}) and (\ref{Df-Df}) that
\[
\begin{split}
&\lVert Df_{n}(x)-Df_{n}(y)  \rVert
\\
&\le \int_n^{n+1}Me^{\bar{\lambda} +\epsilon \lvert r\rvert}\eta e^{-3\epsilon \lvert r \rvert}d\lVert x-y\rVert\, dr \\
&\phantom{\le}+\int_n^{n+1}Me^{\bar{\lambda} +\epsilon \lvert r\rvert}Be^{-4\epsilon \lvert r\rvert} \lVert \phi(r, n;x)-\phi(r, n;y)\rVert \tilde M e^{\epsilon \lvert n\rvert}\, dr \\
&\le \int_n^{n+1}Me^{\bar{\lambda} +\epsilon \lvert r\rvert}\eta e^{-3\epsilon \lvert r \rvert}d\lVert x-y\rVert\, dr \\
&\phantom{\le}+\int_n^{n+1}Me^{\bar{\lambda} +\epsilon \lvert r\rvert}Be^{-4\epsilon \lvert r\rvert} ae^{\epsilon \lvert n\rvert} \tilde M e^{\epsilon \lvert n\rvert}\lVert x-y\rVert \, dr.
\end{split}
\]
This proves the second inequality of (\ref{Dfn0}) holds with $\tilde B:=2ade^{\bar{\lambda}+4\epsilon}B M\tilde M >0$ and the claimed result
(\ref{Dfn0}) is proved.

In what follows, we give a lemma on linearization of $(A_n+f_n)_{n\in \mathbb{Z}}$.
\begin{lemma}\label{IF}
Let $\cA(m, n)$ defined in {\rm (\ref{Amn})} admit a strong exponential dichotomy and let $a_1,...,a_r$ and $b_1,...,b_r$ be given in Lemma~{\rm \ref{ts}} such that {\rm (\ref{DH-cond})} holds. Assume that $\alpha$ is the constant given in the formulation of Theorem~{\rm \ref{thm-mr}} and that $(f_n)_{n\in \Z}$ is a sequence of $C^1$ maps $f_n: \mathbb R^d\to \mathbb R^d$ such that {\rm (\ref{fn0})} and {\rm (\ref{Dfn0})} hold, where $\epsilon \ge 0$ is given in~\eqref{ln2}
and $\tilde \eta>0$ is sufficiently small {\rm (}which tends to $0$ when $\alpha$ tends to its upper bound{\rm )}. Then, there exists a sequence $(h_m)_{m\in \Z}$ of homeomorphisms defined in $\mathbb R^d$
such that
\begin{align}\label{conj_m}
 h_{m+1}\circ (A_m+f_m)=A_m \circ h_m, \quad m\in \Z,
\end{align}
\vspace{-0.5cm}
\begin{align}\label{hnhn}
h_n(x)=x+e^{\epsilon |n|}o(\|x\|^{1+\varrho}),\quad
h_n^{-1}(x)=x+e^{\epsilon |n|}o(\|x\|^{1+\varrho})
\end{align}
for some small $\varrho\in (0,1)$, and that
\begin{align}
&\|h_n(x)-h_n(y)\|\le CLe^{\epsilon \lvert n\rvert} \lVert x-y\rVert^{\alpha}, \label{hnhn-aa1}
\\
&\|h_n^{-1}(x)-h_n^{-1}(y)\|\le CLe^{\epsilon \lvert n\rvert} \lVert x -y\rVert^{\alpha},
\label{hnhn-aa2}
\end{align}
where $L>0$ is constant independent of $\alpha$, for all $x,y\in U_n:=\{u\in\R^d:\|u\|\le C^{-1}e^{-\epsilon|n|}\rho\}$ with a small constant $\rho>0$ independent of $\alpha$.
\end{lemma}
\vspace{-0.3cm}

\noindent Remark that in this lemma if $A_n$ and $f_n$ is independent of $n$ (therefore $\varepsilon=0$ which is given in (\ref{ln2})), then $h_n$ is also independent of $n$ and (\ref{hnhn})-(\ref{hnhn-aa2}) holds with $\varepsilon=0$. The proof of the lemma together with the remark will be postponed to the next section and we continue our proof of Theorem \ref{thm-mr}.
By (\ref{eta}), we understand that $\eta$ (given in the formulation of Theorem \ref{thm-mr}) tends to $0$ when $\alpha$ tends to its upper bound.
Let
\begin{align}\label{def-Htx}
H(t,x):=T(t, n)h_n(\phi(n, t;x)),
\end{align}
for $x\in \R^d$, $t\in [n, n+1)$, $n\in \Z$.
 It is easy to see from \eqref{conj_m} that if $t\mapsto x(t)$ is a solution of~\eqref{943} then  $t\mapsto H(t, x(t))$ is a solution of~\eqref{944}, which proves {\bf (A4)}.
Furthermore, repeating the arguments used to establish~\eqref{912}, we can see that
\begin{equation}\label{5:25}
\|\phi(n, t;x)-\phi(n, t;y)\| \le ae^{\epsilon \lvert n\rvert}\lVert x-y\rVert
\end{equation}
for $t\in [n, n+1)$ and therefore for any
$
x\in V_t:=\{u\in\R^d:\|u\|\le e^{-2\epsilon|t|}\tilde \rho\}
$,
where $\tilde \rho:=(aC)^{-1}e^{-2\epsilon}\rho$ (independent of $\alpha$),
we have
$$
\|\phi(n, t;x)\| \le ae^{\epsilon \lvert n\rvert}\lVert x\rVert
\le C^{-1}e^{-\epsilon|n|}\rho,
$$
implying that $\phi(n, t;x)\in U_n$.
Thus,
using \eqref{nbg} and \eqref{hnhn-aa1}, for $n\in \Z$ and $t\in [n, n+1)$ we have
\begin{align*}
\|H(t,x)-H(t,y)\|&= \| T(t,n)h_n(\phi(n, t;x))-T(t,n)h_n(\phi(n, t;y))\|
\\
&\le \lVert T(t,n)\rVert \cdot \lVert h_n(\phi(n, t;x))-h_n(\phi(n, t;y)) \rVert \\
&\le CMe^{\bar{\lambda}+2\epsilon \lvert n\rvert}\|\phi(n, t;x)-\phi(n, t;y)\|^\alpha
\\
& \le a^\alpha CMe^{\bar{\lambda}+(2\epsilon+\alpha \epsilon) \lvert n\rvert}\|x-y\|^\alpha
\\
&\le a^\alpha CMe^{\bar{\lambda}+2\epsilon+\alpha \epsilon} e^{(2\epsilon+\alpha \epsilon)\lvert t\rvert }\|x-y\|^\alpha
\\
&\le \tilde C e^{(2+\alpha) \epsilon\lvert t\rvert }\|x-y\|^\alpha
\end{align*}
for all $x,y\in V_t$, where $\tilde C:=a^\alpha CMLe^{\bar{\lambda}+2\epsilon+\alpha \epsilon}>0$ is a constant (independent of $\alpha$). This proves the first inequality in {\bf (A2)}.

Moreover,
using~\eqref{nbg}, \eqref{hnhn} and~\eqref{5:25},
for $n\in \Z$ and $t\in [n, n+1)$ we get
\begin{eqnarray}
\lVert H(t,x)-x\rVert
&=&\lVert T(t,n)h_n(\phi(n,t;x))-x\rVert
\nonumber\\
&\le&
\lVert T(t,n)h_n(\phi(n,t;x)) -T(t,n)\phi(n,t;x)\rVert
\nonumber\\
& & 
+\lVert T(t,n)\phi(n,t;x)-T(t,n)T(n,t)x\rVert
\nonumber\\
&\le&
e^{2\epsilon \lvert n\rvert}o(\lVert \phi(n, t;x)\rVert^{1+\varrho})
\nonumber\\
& &
+Me^{\bar{\lambda} +\epsilon \lvert n\rvert}\lVert \phi(n,t;x)-T(n,t)x\rVert
\nonumber\\
&\le&
e^{(3+\varrho)\epsilon \lvert t\rvert} o(\lVert x\rVert^{1+\varrho})
\nonumber\\
& & 
+Me^{\bar{\lambda} +\epsilon \lvert t\rvert}\lVert \phi(n,t;x)-T(n,t)x\rVert.
\label{H-x}
\end{eqnarray}
On the other hand, by \eqref{nbg}, {\bf (F2)}, {\bf (F4)} and~\eqref{5:25} we have
\[
\begin{split}
&\lVert \phi(n,t;x)-T(n,t)x\rVert
\\
&\le \int_n^{n+1}\lVert T(n,s)f(s, \phi(n,s;x))\rVert\, ds
\\
&\le \int_n^{n+1}Me^{\bar{\lambda}+\epsilon \lvert s\rvert}\, \sup_{\theta\in (0,1)}\|D_xf(s, \theta \phi(n,s;x))\|\,\|\phi(n,s;x)\|\,ds
\\
&\le \int_n^{n+1}Me^{\bar{\lambda}+\epsilon \lvert s\rvert}B e^{-4\varepsilon |s|}\|\phi(n,s; x)\|^2 ds
\\
&\le \int_n^{n+1}Me^{\bar{\lambda}+\epsilon \lvert s\rvert}B e^{-4\varepsilon |s|}a^2e^{2\epsilon \lvert n\rvert}\|x\|^2 ds
\\
&\le a^2MBe^{\bar{\lambda}+4\epsilon}\lVert x\rVert^2,
\end{split}
\]
which together with (\ref{H-x}) implies that
\begin{align*}
H(t,x)
&=x+e^{(3+\varrho)\epsilon \lvert t\rvert} o(\lVert x\rVert^{1+\varrho})
+O(\|x\|^2)
\\
&=x+e^{(3+\varrho)\epsilon \lvert t\rvert} o(\lVert x\rVert^{1+\varrho}).
\end{align*}
This proves the first inequality in {\bf (A1)}.

Similarly, we define $G\colon \R \times \R^d \to \R^d$ by
\[
G(t,x)=\phi(t, n;h_n^{-1}(T(n, t)x)),
\]
for $x\in \R^d$, $t\in [n, n+1)$, $n\in \Z$. Again, it is easy to verify that $G$ satisfies the second equality in {\bf (A1)},
the second inequality in {\bf (A2)} and {\bf (A5)}. Finally, we check that
\[
\begin{split}
H(t, G(t, x)) &=T(t, n)h_n(\phi(n, t; G(t,x))) \\
&=T(t,n)h_n (\phi(n, t; \phi(t, n;h_n^{-1}(T(n, t)x)) ) \\
&=T(t,n)h_n (h_n^{-1}(T(n,t)x)) \\
&=T(t,n)T(n, t)x \\
&=x
\end{split}
\]
for each $x\in \R^d$, $t\in [n, n+1)$ and $n\in \Z$. Hence,
\[
H(t, G(t, x))=x \quad \text{for $t\in \R$ and $x\in \R^d$.}
\]
Similarly, one can show that
\[
G(t, H(t, x))=x \quad \text{for $t\in \R$ and $x\in \R^d$.}
\]
This proves {\bf (A3)} and the proof of the theorem is completed.
\qquad$\Box$

A special case of (\ref{943}) is the autonomous system
\begin{equation}\label{AuEq}
x'=Ax+f(x),\quad x\in \mathbb{R}^d,
\end{equation}
where $A$ is a $d\times d$ constant matrix and has $d$ complex eigenvalues $\mu_1,...,\mu_d$ and
$f(0)=0$,
$Df(0)=0$.
One can see easily that in this autonomous case the spectral bound condition (\ref{DH-cond}) holds automatically. Moreover,  (\ref{ED1})-(\ref{nbg}) hold with $\varepsilon=0$. Thus, by Theorem \ref{thm-mr} we obtain the following.

\begin{theorem}\label{thm-aut}
Let the matrix $A$ of system {\rm (\ref{AuEq})} be hyperbolic, i.e.,
 $$
 {\rm Re} ~\mu_1\le \cdots \le {\rm Re} ~\mu_p<0< {\rm Re} ~\mu_{p+1}\le \cdots \le {\rm Re} ~\mu_d,
 $$
where $1\le p\le d-1$ and {\rm Re} denotes the real part of a complex number, and let $f$ be locally $C^{1,1}$ {\rm (}i.e., $f$ is $C^1$ and $Df$ is Lipschitz near the origin{\rm )} such that $f(0)=0$, $Df(0)=0$.
Then, for any $\alpha\in \mathbb{R}$ satisfying
 \begin{align*}
    0<\alpha<\min\Big\{\frac{{\rm Re}~\mu_{p+1}-{\rm Re}~\mu_p}{{\rm Re}~\mu_d},\frac{{\rm Re}~\mu_{p+1}-{\rm Re}~\mu_{p}}{-{\rm Re}~\mu_1}\Big\},
 \end{align*}
 there exist a small neighborhood $V\subset\R^d$ {\rm (}the diameter of $V$ tends to 0 when $\alpha$ tends to its upper bound{\rm )} and a map $\tilde H :V\to \R^d$ such that
 \vspace{-0.2cm}
\begin{itemize}
\item[{\bf(B1)}] $\tilde H(x)=x+o(\lVert x\rVert^{1+\varrho})$, $\tilde H^{-1}(x)=x+o(\lVert x\rVert^{1+\varrho})$ with a small constant $\varrho\in (0,1)$;
\vspace{0.2cm}
\item[{\bf(B2)}] $\|\tilde H(x)-\tilde H(y)\|\le \tilde C \|x-y\|^\alpha$ and $\|\tilde H^{-1}(x)-\tilde H^{-1}(y)\|\le \tilde C\|x-y\|^\alpha$ for all $x,y\in V$, where $\tilde C>0$ is a constant independent of $\alpha$;
\item[{\bf(B3)}] $e^{A t}\tilde H(x)=\tilde H(\phi(t, 0;x))$, where $\phi(t,0;x)$ is the solution $x(t)$ of~\eqref{AuEq} such that $x(0)=x$, i.e., $\tilde H$ is a conjugacy between {\rm (\ref{AuEq})} and the linear system $x'=Ax$.
\end{itemize}
\vspace{-0.4cm}
\end{theorem}

\noindent We remark  that the difference between the main result in~\cite{CDS}  and ours is the most  obvious in the above described autonomous case. Indeed,  while~\cite{CDS} still needs the non-resonant conditions
up to order $k$,  we do not need any non-resonant conditions in Theorem~\ref{thm-aut}. Moreover, we stress that Theorem~\ref{thm-aut}  is the first result that gives a rigorous proof for simultaneously differentiable and H\"older linearization of hyperbolic systems without any non-resonant conditions because van Strien's proof~\cite{Stri-JDE90} was pointed out to be wrong, as we already mentioned in the introduction. Notice that this theorem is a continuous-time version of van Strien's result. In the discrete-time case, the remark given just below Lemma \ref{IF} shows that Lemma \ref{IF} with $A_n$ and $f_n$ independent of $n$ is has the same framework as van Strien's result.

{\it Proof of Theorem {\rm \ref{thm-aut}}}. For any given small constant $\eta>0$, there is a small neighborhood $V_\eta\subset \mathbb{R}^d$ of the origin such that one can use a smooth cut-off function defined in $\R^d$ (i.e., a smooth function which is equal to 1 in $V_\eta$ and is equal to 0 outside a neighborhood of $V_\eta$) to extend the locally defined $C^{1,1}$ map $f$ to a global one satisfying
\begin{align*}
\|D f(x)\| \le \eta\quad{\rm and}\quad\|Df(x)-Df(y)\|\le B\|x-y\|,\quad \forall x,y\in \mathbb{R}^n,
\end{align*}
where $B>0$ is a constant (see e.g. \cite{ZZJ}). Notice that the diameter of $V_\eta$ tends to 0 when $\eta$ tends to 0. On the other hand,
by Theorem \ref{thm-mr} we see that $\eta$ tends to $0$ when $\alpha$ tends to its upper bound.
Hence, one concludes that the diameter of $V_\eta$ tends to 0 when $\alpha$ tends to its upper bound.

Next, one checks that {\bf (F1)}-{\bf (F4)} hold with $f(t,x)$ replacing by $f(x)$ and with $\varepsilon=0$. Then by Theorem \ref{thm-mr} and (\ref{def-Htx}) we obtain
\[
H(t,x)=T(t, n)h_n(\phi(n, t;x))
\]
for $x\in \R^d$, $t\in [n, n+1)$, $n\in \Z$, which satisfies {\bf (A1)}-{\bf (A2)} with $\varepsilon=0$. Notice that in the autonomous case
\begin{align}\label{TATA}
T(t,n)=e^{A(t-n)},\quad \phi(n, t;x)=\phi(n-t, 0;x),
\end{align}
and therefore
\begin{align}\label{THTH}
H(t,x)=T(t, n)h_n(\phi(n, t;x))
=e^{A(t-n)}h(\phi(n-t, 0;x)).
\end{align}
Then, $\tilde H:\R^n\to \R^n$ can be defined by
\begin{align*}
\tilde H(x):=\int_{0}^{1}e^{As}h(\phi(-s, 0;x))ds
=\int_{n}^{n+1}H(t,x)dt.
\end{align*}
Similarly, $\tilde H^{-1}(x)$ can be obtained by $G(t,x)$ and
by {\bf (A1)}-{\bf (A2)} one verifies {\bf (B1)}-{\bf (B2)} for $V:=V_\eta\cap \{u\in\R^d:\|u\|\le \tilde \rho\}$ with a small constant $\tilde \rho>0$ given in Theorem~\ref{thm-mr}. It is clear that the diameter of $V$ tends to 0 when $\alpha$ tends to its upper bound by the last sentence of the previous paragraph.

Moreover, it also follows from (\ref{TATA}) that $A_n=e^A$, as seen in (\ref{an}), and $f_n=\phi(1,0;\cdot)-e^A$,
as seen in (\ref{xfn}), both of which are independent of $n$.
Thus, $h_n$ obtained in Lemma \ref{IF} can be independent of $n$ by the remark given just below Lemma \ref{IF}.
This enables us to put $h:=h_n$ and rewrite (\ref{conj_m}) as
\begin{equation}\label{conAh}
 h\circ \phi(1,0;\cdot)=e^A \circ h.
\end{equation}
Then we see that
\begin{align*}
&e^{At}\tilde H(x)
\\
&=\int_{0}^{1}e^{A(s+t)}h(\phi(-s-t, 0; \phi(t, 0;x)))ds
\\
&=\int_{t}^{0}e^{As}h(\phi(-s, 0;\phi(t, 0;x)))ds+\int_{0}^{1+t} e^{As}h(\phi(-s, 0;\phi(t, 0;x)))ds
\\
&=\int_{t}^{0}\!\!e^{A(s+1)}h(\phi(-s\!-\!1, 0;\phi(t, 0;x)))ds
+\!\!\int_{0}^{1+t}\!\!e^{As}h(\phi(-s, 0;\phi(t, 0;x)))ds
\\
&=\int_{1+t}^{1}e^{As}h(\phi(-s, 0;\phi(t, 0;x)))ds +\int_{0}^{1+t}e^{As}h(\phi(-s, 0;\phi(t, 0;x)))ds
\\
&=\int_{0}^{1}e^{As}h(\phi(-s, 0;\phi(t, 0;x)))ds \\
&=\tilde H(\phi(t, 0;x)),
\end{align*}
where we have used that
\begin{align*}
\int_{t}^{0}e^{As}h(\phi(-s, 0;x))ds
&=\int_{t}^{0}e^{A(s+1)}e^{-A}h(\phi(1,0;\phi(-s-1, 0;x)))ds
\\
&=
\int_{t}^{0}e^{A(s+1)}h(\phi(-s-1, 0;x))ds
\end{align*}
since $e^{-A}h(\phi(1,0;y))=h(y)$, as seen from (\ref{conAh}). This proves {\bf (B3)} and the proof of the theorem is completed. \qquad$\Box$


\section{Proof of Lemma \ref{IF}}

{\it Proof of Lemma} \ref{IF}.
Define  a map $F\colon Y_\infty \to Y_\infty$ by
\begin{align}\label{DFx}
 (F(\mathbf x))_n:=A_{n-1}x_{n-1}+f_{n-1}(x_{n-1}), \quad \mathbf x=(x_n)_{n\in \Z} \in Y_\infty.
\end{align}
By the same argument as in \cite[Claims 3 and 4]{DZZ}, we can see that
\begin{itemize}
\item $F$ is well-defined and differentiable such that
\[
DF(\mathbf x) \xi=(A_{n-1}\xi_{n-1}+Df_{n-1}(x_{n-1})\,\xi_{n-1})_{n\in \Z}
\]
for each $\mathbf x=(x_n)_{n\in \Z}$ and $\mathbf \xi=(\xi_n)_{n\in \Z} \in Y_\infty$;
\vspace{0.2cm}
\item $F$ is $C^{1, 1}$, which means that
 \[
  \sup_{\mathbf x\neq \mathbf y} \frac{\lVert  DF(\mathbf x)-DF(\mathbf y) \rVert}{\lVert \mathbf x-\mathbf y\rVert}\le C\tilde B<\infty;
 \]
\item $\lVert DF(\mathbf x)-\mathbb A\rVert \le C\tilde \eta$ for all $\mathbf x\in Y_\infty$.
\end{itemize}
Hence, $\mathbf 0:=(0)_{n\in \Z}$ is a hyperbolic fixed point of~$F$ since one sees from (\ref{DFx}) that $DF(\mathbf 0) =\mathbb A$ and $\mathbb A$ is hyperbolic, i.e. $\sigma (\mathbb A)\cap S^1=\emptyset$. Then we have the following lemma on smooth linearization of $F$.
\begin{lemma}\label{lm-F}
Let $F$ and $\mathbb{A}$ be given above and assume that the numbers $a_i$ and $b_i$, given in the statement of Lemma~{\rm \ref{ts}}, satisfy {\rm (\ref{DH-cond})}. Then, for the constant $\alpha$ given in the formulation of Theorem~{\rm \ref{thm-mr}}, there exists a homeomorphism $\Phi: X\to X$ such that
\begin{equation}\label{conj}
 \Phi \circ F=\mathbb A \circ \Phi,
\end{equation}
where $\Phi$ and $\Phi^{-1}$ satisfy that
\begin{align}
\Phi({\bf x})={\bf x}+O(\|{\bf x}\|^{1+\varrho}),\quad \Phi^{-1}({\bf x})={\bf x}+O(\|{\bf x}\|^{1+\varrho})\quad {\rm as}~ \|{\bf x}\|\to 0,
\label{bbb}
\end{align}
for some small $\varrho\in (0,1)$ and are both $\alpha$-H\"older continuous in $\{{\bf x}\in X: \|{\bf x}\|\le\rho\}$ with a small constant $\rho>0$ independent of $\alpha$, i.e.,
\begin{align*}
\|\Phi({\bf x})-\Phi({\bf y})\|\le L\|{\bf x}-{\bf y}\|^{\alpha},\quad
\|\Phi^{-1}({\bf x})-\Phi^{-1}({\bf y})\|\le L\|{\bf x}-{\bf y}\|^{\alpha},
\end{align*}
where $L>0$ is a constant independent of $\alpha$.
\end{lemma}

{\it Proof of Lemma} \ref{lm-F}. Since $\sigma (\mathbb A)\cap S^1=\emptyset$ as mentioned before, the space $Y_\infty$ has a direct decomposition
$$
Y_\infty:=Y_s\oplus Y_u,
$$
where $Y_s$ and $Y_u$ correspond to the spectra
$$
\bigcup_{i=1}^k \{z\in \mathbb{C}: a_i\le|z|\le b_i\} \quad{\rm and}~~~~~
\bigcup_{i=k+1}^r \{z\in \mathbb{C}: a_i\le|z|\le b_i\},
$$
respectively. Thus ${\bf x}={\bf x}_s+{\bf x}_u$ where ${\bf x}_s\in Y_s$ and
${\bf x}_u\in Y_u$. Let $\pi_s$ and $\pi_u$ be projections such that
$$
\pi_s {\bf x}:={\bf x}_s\quad{\rm and}\quad\pi_u {\bf x}:={\bf x}_u
$$
and let $\|{\bf x}\|=\|\pi_s{\bf x}\|+\|\pi_u{\bf x}\|$. Denote $\mathbb{A}_s:=\mathbb{A}|_{Y_s}$, $\mathbb{A}_u:=\mathbb{A}|_{Y_u}$ and $\tilde f:=F-\mathbb{A}$.
By the discussion given in the proof of \cite[Theorem 1]{ZZJ}, we understand that the key step of the proof is to solve the functional equation
\begin{align}
&q_n({\bf x},\xi_s)=\mathbb{A}_s^n(\xi_s-\pi_s
{\bf x})
\nonumber\\
&\quad+\sum_{i=0}^{n-1}\mathbb{A}_s^{n-i-1}\big\{\pi_s
{\tilde f}(q_i({\bf x},\xi_s)+F^i({\bf x}))-\pi_s {\tilde f}(F^i({\bf x}))\big\}
\nonumber\\
&\quad-\sum_{i=n}^{+\infty}\mathbb{A}_u^{n-i-1}\big\{\pi_u
{\tilde f}(q_i({\bf x},\xi_s)+F^i({\bf x}))-\pi_u {\tilde f}(F^i({\bf x}))\big\},
\quad \forall n\ge 0,
\label{eqns-foli-seq}
\end{align}
with $q_n:Y_\infty\times Y_s\to Y_\infty$ unknown, which can be used to define the stable foliation of $Y_\infty$ under $F$. Once we find a smooth solution $(q_n)_{n\ge 0}$ of Eq. (\ref{eqns-foli-seq}), the corresponding stable foliation that has the same smoothness as $q_0$ can be constructed. Notice that an unstable invariant foliation can be obtained by considering the stable one of the inverse $F^{-1}$. Then, using a transformation that has the same smoothness as the stable and unstable foliations,
we may decouple $F$ into a $C^{1,1}$ contraction and a $C^{1,1}$ expansion.
Finally, smooth linearization theorem for contractions can be use to complete the proof of this lemma.

Following the above strategy, in order to solve Eq. (\ref{eqns-foli-seq}) we know from
\cite[Theorems 2.1-2.2]{ChHaTan-JDE97} that equation (\ref{eqns-foli-seq}) has a unique $C^0$ solution $(q^*_n)_{n\ge 0}$
such that
$
\sup_{n\ge 0}\{r^{-n}\|q^*_n({\bf x},\xi_s)\|\}<\infty
$
for any constant $r\in (b_k,a_{k+1})$ and for every fixed $({\bf x},\xi_s)\in Y_\infty\times Y_s$. Then
\cite[Lemma 7.1]{ZLZ} tells that $q_0^*$ satisfies
\begin{align}
\sup_{({\bf x},\xi_s)\in\Omega\backslash\{(0,0)\}}\frac{\|q^*_0({\bf x},\xi_s)-(\xi_s-\pi_s
{\bf x})\|}{\|({\bf x},\xi_s)\|^{1+\varrho}}<\infty
\label{q0LL22}
\end{align}
for a small constant $\varrho\in (0,1)$, where $\Omega\subset Y_\infty\times Y_s$ is a small neighborhood of the origin $(0,0)$ in the space $Y_\infty\times Y_s$.
In what follows, we further show that $q_0^*$ is H\"older continuous.
In fact, since $(q^*_n)_{n\ge 0}$ is a solution of equation (\ref{eqns-foli-seq}), we have
\begin{align}
q^*_n({\bf x},\xi_s)=~&\mathbb{A}_s^n(\xi_s-\pi_s
{\bf x})
\nonumber\\
&+\sum_{i=0}^{n-1}\mathbb{A}_s^{n-i-1}\big\{\pi_s
{\tilde f}(q^*_i({\bf x},\xi_s)+F^i({\bf x}))-\pi_s {\tilde f}(F^i({\bf x}))\big\}
\nonumber\\
&-\sum_{i=n}^{+\infty}\mathbb{A}_u^{n-i-1}\big\{\pi_u
{\tilde f}(q^*_i({\bf x},\xi_s)+F^i({\bf x}))-\pi_u {\tilde f}(F^i({\bf x}))\big\}
\label{eqns-foli}
\end{align}
for $n\ge 0$. Choose constants $\lambda_s^+,\lambda_u^-,\gamma_s,\gamma_u\in (b_k,a_{k+1})$ and $\lambda_u^+\in (b_r,\infty)$ such that
$$
b_k<\lambda_s^+<\gamma_s<1<\gamma_u<\lambda_u^-<a_{k+1}\quad
{\rm and}\quad \gamma_s\gamma_u^{-1}(\lambda_u^+)^{\alpha}<1,
$$
the second of which is possible due to (\ref{def-a}). By \cite[Theorem 5]{R-S-JDE04}, one can choose appropriate equivalent norms
in $Y_\infty$ such that
\begin{align*}
\|\mathbb{A}_s\|< \lambda_s^+, \quad
\|\mathbb{A}_u^{-1}\|< 1/\lambda_u^-, \quad \|\mathbb{A}\|=\|\mathbb{A}_u\|< \lambda_u^+.
\end{align*}
Thus,
\begin{align*}
\gamma_s^{-n}\|q^*_n({\bf x},\xi_s)\|
&\le \gamma_s^{-n}\|\mathbb{A}_s\|^n(\|\xi_s\|+\|\pi_s
{\bf x}\|)
\\
&\qquad+\gamma_s^{-1}\sum_{i=0}^{n-1}\gamma_s^{-(n-i-1)}\|\mathbb{A}_s\|^{n-i-1}
\\
&\qquad\quad\cdot\gamma_s^{-i}\big\|\pi_s
{\tilde f}(q^*_i({\bf x},\xi_s)+F^i({\bf x}))-\pi_s {\tilde f}(F^i({\bf x}))\big\|
\\
&\qquad+\gamma_s^{-1}\sum_{i=n}^{+\infty}\gamma_s^{-(n-i-1)}\|\mathbb{A}_u^{-1}\|^{-(n-i-1)}
\\
&\qquad\quad\cdot
\gamma_s^{-i}\big\|\pi_u
{\tilde f}(q^*_i({\bf x},\xi_s)+F^i({\bf x}))-\pi_u {\tilde f}(F^i({\bf x}))\big\|
\\
&\le\Big(\frac{\lambda_s^+}{\gamma_s}\Big)^n(\|{\bf x}\|+\|\xi_s\|)
+\gamma_s^{-1}\sum_{j=0}^{\infty}\Big\{\Big(\frac{\lambda_s^+}{\gamma_s}\Big)^{j}
+\Big(\frac{\gamma_s}{\lambda_u^-}\Big)^{j}\Big\}
\\
&\qquad \cdot \sup_{{\bf z}\in Y_\infty}\|Df({\bf z})\|\sup_{i\ge 0}\{\gamma_s^{-i}\|q^*_i({\bf x},\xi_s)\|\}
\\
&\le\|{\bf x}\|+\|\xi_s\|+CK\tilde \eta\sup_{i\ge 0}\{\gamma_s^{-i}\|q^*_i({\bf x},\xi_s)\|\}
\end{align*}
because $\|Df({\bf z})\|\le C\tilde \eta$,
where $\tilde \eta>0$ is small enough such that $CK\tilde \eta<1/4$.
It implies that for any small $\delta>0$, we can choose $\Omega:=\{({\bf x},\xi_s)\in Y_\infty\times Y_s: \|{\bf x}\|\le \delta/4, ~\|\xi_s\|\le \delta/4\}$ such that
\begin{align*}
\sup_{n\ge 0}\{\gamma_s^{-n}\|q^*_n({\bf x},\xi_s)\|\}\le 2(\|{\bf x}\|+\|\xi_s\|)\le \delta,
\qquad \forall ({\bf x},\xi_s)\in \Omega.
\end{align*}
Let $[{\tilde f}(\cdot)]_{\bf y}^{\bf x}:={\tilde f}({\bf x})-{\tilde f}({\bf y})$. Then for all $({\bf x},\xi_s),({\bf y},\xi_s)\in \Omega$
\begin{align}
&\gamma_u^{-n}\|q^*_n({\bf x},\xi_s)-q^*_n({\bf y},\xi_s)\|
\nonumber\\
&\le \Big(\frac{\lambda_s^+}{\gamma_u}\Big)^n\|{\bf x}-{\bf y}\|+\gamma_u^{-1}
  \sum_{j=0}^{\infty}\Big\{\Big(\frac{\lambda_s^+}{\gamma_u}\Big)^{j}
  +\Big(\frac{\gamma_u}{\lambda_u^-}\Big)^{j}\Big\}
  \nonumber\\
  &~~\quad
  \cdot \sup_{i\ge 0}\Big\{\gamma_u^{-i}\Big\|\Big[\int_0^1 D{\tilde f}(tq^*_i(\cdot,\xi_s)+F^i(\cdot))q^*_i(\cdot,\xi_s)dt\Big]_{{\bf y}}^{\bf x}\Big\|\Big\}
\nonumber\\
&\le \|{\bf x}-{\bf y}\|+K\sup_{i\ge 0}\Big\{\gamma_u^{-i}\sup_{t\in
  (0,1)}\Big\|\Big[D{\tilde f}(tq^*_i(\cdot,\xi_s)+F^i(\cdot))\Big]_{{\bf y}}^{\bf x}\Big\|\,\|q^*_i({\bf x},\xi_s)\|\Big\}
  \nonumber\\
  &~~\quad +K\sup_{i\ge 0}\Big\{\gamma_u^{-i}\sup_{t\in (0,1)}\|D{\tilde f}(tq^*_i({\bf y},\xi_s)+F^i({\bf y}))\|\,\|q^*_i({\bf x},\xi_s)-q^*_i({\bf y},\xi_s)\|\Big\}
\nonumber\\
&\le \|{\bf x}-{\bf y}\|+K\delta\sup_{i\ge 0}\sup_{t\in
  (0,1)}\bigg\{\gamma_s^{i}\gamma_u^{-i}\Big(\Big\|\Big[D{\tilde f}(tq^*_i(\cdot,\xi_s)
  +F^i({\bf x}))\Big]_{{\bf y}}^{\bf x}\Big\|
  \nonumber\\
  &~~\quad\!+\!\Big\|\Big[D{\tilde f}(tq^*_i({\bf y},\xi_s)+F^i(\cdot))\Big]_{{\bf y}}^{\bf x}\Big\|\Big)\bigg\}
  \!+\!CK\tilde \eta  \sup_{i\ge 0}\{\gamma_u^{-i}\|q^*_i({\bf x},\xi_s)-q^*_i({\bf y},\xi_s)\|\}
\nonumber\\
&\le \|{\bf x}-{\bf y}\|+KC\tilde B\delta\sup_{i\ge 0}\bigg\{\gamma_s^{i}\gamma_u^{-i}
  \Big(\|q^*_i({\bf x},\xi_s)-q^*_i({\bf y},\xi_s)\|+(2C\tilde \eta)^{1-\alpha}
  \nonumber\\
  &~~\quad \cdot (C\tilde B)^{\alpha}\|F^i({\bf x})-F^i({\bf y})\|^{\alpha}\Big)\bigg\}
  + CK\tilde \eta  \sup_{i\ge 0}\{\gamma_u^{-i}\|q^*_i({\bf x},\xi_s)-q^*_i({\bf y},\xi_s)\|\}
  \nonumber\\
&\le \|{\bf x}-{\bf y}\|^\alpha+\frac{1}{4}\sup_{i\ge 0}\Big\{\gamma_u^{-i}
  \|q^*_i({\bf x},\xi_s)-q^*_i({\bf y},\xi_s)\|
\nonumber\\
&~~\quad +\gamma_s^{i}\gamma_u^{-i}(\lambda_u^++C\tilde \eta)^{i\alpha}\|{\bf x}-{\bf y}\|^{\alpha}\Big\}+ \frac{1}{4}  \sup_{i\ge 0}\{\gamma_u^{-i}\|q^*_i({\bf x},\xi_s)-q^*_i({\bf y},\xi_s)\|\}
\nonumber\\
&\le \frac{5}{4}\|{\bf x}-{\bf y}\|^{\alpha}+\frac{1}{2}\sup_{i\ge 0}\{\gamma_u^{-i}\|q^*_i({\bf x},\xi_s)-q^*_i({\bf y},\xi_s)\|\},
\label{long-f}
\end{align}
where $\delta,\tilde \eta>0$ are small enough such that
$$
K(C\tilde B)^2\delta<1/4,\quad CK\tilde \eta<1/4,\quad \gamma_s\gamma_u^{-1}(\lambda_u^++C\tilde \eta)^{\alpha}<1.
$$
Notice that the constants $K$ (depends only on the linear part $\mathbb A$), $\tilde B$ (the Lipschitz constant of $D f_n$) and $C$ (see (\ref{ln2})) are all independent of $\alpha$ and so does the constant $\delta$. However, $\tilde \eta$ depends on $\alpha$ since we need $\gamma_s\gamma_u^{-1}(\lambda_u^++C\tilde \eta)^{\alpha}<1$ and therefore $\tilde \eta$ tends to $0$ when $\alpha$ tends to its upper bound.

It follows from (\ref{long-f}) that
$
\sup_{i\ge 0}\{\gamma_u^{-i}\|q^*_i({\bf x},\xi_s)-q^*_i({\bf y},\xi_s)\|\}\le 3\|{\bf x}-{\bf y}\|^{\alpha}
$
and therefore
\begin{align}\label{HOLD}
\|q^*_0({\bf x},\xi_s)-q^*_0({\bf y},\xi_s)\|\le 3\|{\bf x}-{\bf y}\|^{\alpha}.
\end{align}
The locally $\alpha$-H\"older continuity of $q^*_0({\bf x},\xi_s)$ in $\xi_s$ is clear, i.e.,
\begin{align}\label{HOLD2}
\|q^*_0({\bf x},\xi_s)-q^*_0({\bf x},\tilde \xi_s)\|\le L\|\xi_s-\tilde \xi_s\|^{\alpha}
\end{align}
for a constant $L>0$ since it is actually $C^1$ in $\xi_s$ by \cite[Theorem 1.1]{ChHaTan-JDE97}.
Hence, in view of (\ref{q0LL22}), (\ref{HOLD}), (\ref{HOLD2}) and the discussion given in the proof of \cite[Theorem 1]{ZZJ}, we understand that Theorem 3.1 given
in \cite{Tan-JDE00} can be applied to find a
homeomorphism $\Psi:Y_\infty\to Y_\infty$, which satisfies that
\begin{align*}
\Psi({\bf x})={\bf x}+O(\|{\bf x}\|^{1+\varrho}),\quad \Psi^{-1}({\bf x})={\bf x}+O(\|{\bf x}\|^{1+\varrho})\quad {\rm as}~ \|{\bf x}\|\to 0,
\end{align*}
and that both $\Psi$ and $\Psi^{-1}$ are $\alpha$-H\"older continuous in $\{{\bf x}\in X: \|{\bf x}\|\le\rho\}$ with $\rho:=\delta/4$ (independent of $\alpha$), such that the equality
\begin{align*}
\Psi\circ F&=F_s\circ \pi_s\Psi+F_u\circ \pi_u\Psi
\end{align*}
holds. Here, the maps $F_s: Y_s\to Y_s$ and $F_u:Y_u\to Y_u$ are defined by
\begin{align*}
F_s:=\pi_sF\circ({\rm id}_s+g_s),
~~~~~
F_u:=\pi_uF\circ({\rm id}_u+g_u),
\end{align*}
where ${\rm id}_j$'s are identity mappings in $Y_j$'s for $j=s,u$, and the graphs of $g_j:Y_j\to Y_j$  for $j=s,u$ are $C^{1,1}$ stable and unstable invariant manifolds, respectively. Therefore,
one checks that $F_j$'s are $C^{1,1}$ maps such that $DF_s(0)=\mathbb{A}|_{Y_s}$ and
$DF_u(0)=\mathbb{A}|_{Y_u}$.
 Then, by \cite[Lemma 10]{ZZJ} we see that under the spectral bound condition (\ref{DH-cond}) there exist neighborhoods $U_s\subset Y_s$ and
$U_u\subset Y_u$ of the origins and homeomorphisms
$\psi_s:U_s\to Y_s$ and $\psi_u: U_u\to Y_u$, both of which together with their inverses
are $C^{1,\beta}$ with a small constant $\beta\in (0,1)$ such that
\begin{eqnarray*}
\psi_s\circ F_s=\mathbb{A}_s\circ \psi_s, \quad \psi_u\circ
F_u=\mathbb{A}_u\circ \psi_u.
\end{eqnarray*}
This enables us to define a homeomorphism $\Phi$ by
\begin{align*}
\Phi:=(\psi_s\circ \pi_s+\psi_u\circ \pi_u)\circ \Psi,
\end{align*}
which satisfies that $\Phi\circ F=\mathbb{A}\circ \Phi$.
One can further check that $\Phi^{-1}=\Psi^{-1}\circ (\psi_s^{-1}\circ
\pi_s+\psi_u^{-1}\circ \pi_u)$ and that both $\Phi$ and $\Phi^{-1}$ are $\alpha$-H\"older continuous and satisfy (\ref{bbb}).
The proof of Lemma \ref{lm-F} is completed. $\Box$

We continue to prove Lemma \ref{IF}.
For a fixed $n\in \Z$ and $v\in \mathbb R^d$, define $\mathbf x^n=(x_m)_{m\in \Z}$ by $x_n=x$ and $x_m=0$ for $m\neq n$.
Let $h_n(x):=(\Phi(\mathbf x^n))_n$. It follows readily from~\eqref{conj} that~\eqref{conj_m} holds. Furthermore, we see that
\begin{align*}
   \frac{\lVert h_n(x)-x\rVert}{\lVert x\rVert^{1+\varrho}} \displaybreak[0]
   \le Ce^{\epsilon \lvert n\rvert}\, \frac{\lVert h_n(x)
-x\rVert_n}{ \lVert x\rVert_n^{1+\varrho}} \displaybreak[0]
   \le  Ce^{\epsilon \lvert n\rvert} \,\frac{\lVert \Phi(\mathbf x^n)
-\mathbf x^n \rVert}{\lVert \mathbf x^n\rVert^{1+\varrho}}.
\end{align*}
Letting $\lVert x\rVert \to 0$, we have $\lVert \mathbf x^n\rVert \to 0$ and therefore for every $n$
$$
 \frac{\lVert h_n(x)-x\rVert}{Ce^{\epsilon \lvert n\rvert}\lVert x\rVert^{1+\varrho}}\to 0
$$
by (\ref{bbb}), which proves the first equality of (\ref{hnhn}).
The $\alpha$-H\"older smoothness of $h_n$ can be implied by the $\alpha$-H\"older smoothness
of $\Phi$ in Lemma \ref{lm-F}, where $\alpha$ is given in (\ref{def-a}). In fact, from (\ref{ln2})
we understand that if $\|x\|\le C^{-1}e^{-\epsilon|n|}\rho$ then $\|\mathbf x^n\|\le \rho$ with small constant $\rho>0$. Therefore for any $x,y\in U_n$, which is defined in the formulation of Lemma \ref{IF}, we see that
\begin{equation*}
 \begin{split}
  \lVert h_n(x)-h_n(y)\rVert
  &\le \lVert h_n(x)-h_n(y)\rVert_n
  \le \lVert \Phi(\mathbf x^n) -\Phi(\mathbf y^n)\rVert
  \\
  &\le L\lVert \mathbf x^n -\mathbf y^n\rVert^{\alpha}=L\lVert x -y\rVert_n^{\alpha}
  \\
  &\le CLe^{\epsilon \lvert n\rvert} \lVert x -y\rVert^{\alpha},
 \end{split}
\end{equation*}
which proves (\ref{hnhn-aa1}).
Furthermore, we see that
\[
 h_n^{-1}(v)=(\Phi^{-1}(\mathbf v^n))_n \quad \text{for $v\in \mathbb R^d$ and $n\in \Z$.}
\]
Hence, one can repeat the above arguments and show that $h_n^{-1}$ satisfies the second equality of (\ref{hnhn}) and (\ref{hnhn-aa2}).
The proof of the lemma is completed. \qquad $\Box$

The remark given just below Lemma \ref{IF} can be seen easily by Lemma \ref{lm-F}, which shows that if the given system is independent of $n$ then so is the conjugacy.

\section{Infinite-dimensional case}

In this section we briefly discuss how one can extend our results to the case of infinite dimension under suitable additional  assumptions. Let $X$ be an arbitrary Banach space and denote by $B(X)$ the space of all bounded operators on $X$. We now consider equations~\eqref{943} and~\eqref{944}, where
$A\colon \R \to B(X)$ is a continuous map and $f\colon \R \times X\to X$  is a continuous map such that $f(t, \cdot)\colon X\to X$ is $C^1$ for each $t\in \R$. Finally, let $T(t,s)$ denote the evolution family corresponding to~\eqref{944}.

Let us now assume that~\eqref{944} admits a nonuniform strong exponential dichotomy. This means that there exist projections $P(t)$, $t\in \R$ on $X$ such that (4), (5) and~\eqref{nbg} hold with some $M, \lambda, \bar{\lambda}>0$, $\lambda \le \bar{\lambda}$ and $\epsilon \ge 0$. Now one can construct  the family of
$\lVert \cdot \rVert_t$, $t\in \R$ on $X$ as in Section~\ref{sec.2}. Moreover, set
\[
A_n=T(n+1, n), \quad n\in \Z,
\]
and consider a bounded linear operator $\mathbb A\colon Y_\infty \to Y_\infty$ defined by~\eqref{mathbbA} on
\[
Y_\infty:=\bigg{\{} \mathbf x=(x_n)_{n\in \Z} \subset X: \sup_{n\in \Z}\lVert x_n\rVert_n <\infty \bigg{\}}.
\]
Furthermore, suppose that $\sigma (\mathbb A)$ is given by~\eqref{nxx}. One can now repeat all of our  previous arguments and establish the version of Theorem~\ref{thm-mr} in this new setting. We stress that the statement and the proof remain unchanged.

{\bf Acknowledgement:}

The authors are ranked in alphabetic order of their names and should be treated equally.
The author Davor Dragi\v cevi\' c is supported in part by the Croatian Science Foundation under the project IP-2014-09-2285
and by the University of Rijeka under the project number uniri-prirod-18-9.
The author Weinian Zhang is supported by
NSFC grants \# 11771307 and \# 11831012.
The author Wenmeng Zhang is supported by NSFC grants \# 11671061 and NSF-CQ \# cstc2018jcyjAX0418.

We also would like to express our sincere appreciation to the referees for their great patience and valuable suggestions and questions.


\bibliographystyle{amsplain}

\end{document}